\documentclass[11pt]{amsart}
\usepackage{amsmath,amssymb,amsthm}
\usepackage[latin1]{inputenc}
\usepackage{version,tabularx,multicol}
\usepackage{graphicx,float,psfrag}
\usepackage{stmaryrd}

\theoremstyle{plain}
\newtheorem{theo}{Theorem}[section]

\newtheorem{lem}[theo]{Lemma}
 \newtheorem{lemma}[theo]{Lemma}

\theoremstyle{remark}

 \newtheorem{remark}[theo]{Remark}

\newcommand{\mE}{{\mathbb E}}

\newcommand{\mN}{{\mathbb N}}

\newcommand{\mP}{{\mathbb P}}

\newcommand{\mR}{{\mathbb R}}

 \def\beqlb{\begin{eqnarray}}\def\eeqlb{\end{eqnarray}}
 \def\beqnn{\begin{eqnarray*}}\def\eeqnn{\end{eqnarray*}}

 \def\ar{\!\!&}

 \def\mbb{\mathbb}

 \def\qed{\hfill$\Box$\medskip}

\newcommand{\bcen}{\begin{center}}
\newcommand{\ecen}{\end{center}}
\newcommand{\bgeqn}{\begin{equation}}
\newcommand{\edeqn}{\end{equation}}

\def\mN{{\mbb  N}}

\def\mE{{\mbb  E}}

\def\mP{{\mbb  P}}
\def\mR{{\mbb  R}}

\def\mT{{\mbb  T}}
\def\t{{\mathbf t}}
\def\s{{\mathbf s}}

\def\l{\left}
\def\r{\right}

\newcommand{\cal}{\mathcal}
 \def\ar{\!\!\!&}

\begin{document}

\title[Conditioning trees on maximal out-degree]{Conditioning Galton-Watson trees on
large maximal out-degree}
\date{\today}

\author{Xin He}

\address{Xin He, School of Mathematical Sciences, Beijing Normal University, Beijing 100875, P.R.CHINA}

\email{hexin@bnu.edu.cn}

\begin{abstract}
We propose a new way to condition random trees, that is,
condition random trees to have large maximal out-degree.
Under this new conditioning, we show that conditioned critical Galton-Watson trees converge locally
to size-biased trees with a unique infinite spine.
For the sub-critical case, we obtain local convergence to size-biased trees with a unique infinite node.
We also study tail of the maximal out-degree of sub-critical Galton-Watson trees,
which is essential for the proof of the local convergence.
\end{abstract}



\maketitle

\section{Introduction}\label{sec:introduction}

In the seminal
paper \cite{K86}, Kesten studied the local limit of a critical or sub-critical Galton-Watson tree (GW tree) conditioned to have large height. This limit is the size-biased tree with a unique infinite spine, which we call $\emph{Kesten's tree}$ throughout the present paper. Since then, other conditionings have also been considered: large total
progeny, and large number of leaves. In particular, Jonnsson and Stef\'{a}nsson \cite{JS11} noticed that some sub-critical GW trees conditioned to have large total progeny do not converge to Kesten's tree. They proved that instead those large conditioned trees converge to size-biased trees with a unique infinite node, which we call $\emph{condensation tree}$. Janson \cite{J12} completed this result by proving that any sub-critical GW tree conditioned to have large total progeny converges to the condensation tree, if not to Kesten's tree. In \cite{AD14a,AD14b}, Abraham and Delmas provided a convenient framework to study local limits of conditioned GW trees, then they used this framework to prove essentially all previous results and some new ones. Specifically, in \cite{AD14a} they provided a criterion for local convergence of finite random trees to Kesten's tree, then gave short and elementary proofs of essentially all previous related results and some new ones. In \cite{AD14b}, they provided a criterion for local convergence of finite random trees to condensation tree, then generalized and completed essentially all previous related results by conditioning GW
trees to have a large number of individuals with out-degree in a given set.

However we still feel that there is something missing in this seemingly complete picture. Intuitively speaking, it is natural that when we condition a GW tree to have large height, we get Kesten's tree in the limit, which has infinite height. But the other conditioning, conditioning a GW tree to have large total progeny is not the most natural way to get condensation tree in the limit, and the results do show: Some conditioned sub-critical GW trees converge locally to condensation tree, while many other sub-critical ones to Kesten's tree. Intuitively, it is clear that a more natural conditioning to get condensation tree in the limit, should be conditioning GW trees to have large $\emph{maximal out-degree}$, which is the maximal number of offsprings of all the individuals in the tree.

In the present paper, we condition GW trees to have large maximal out-degree and use the framework in \cite{AD14a,AD14b} to study local limits of large conditioned trees.
We say that a probability distribution $p=(p_0,p_1,p_2,\ldots)$ on non-negative integers is $\emph{bounded}$, if the set $\{n; p_n > 0\}$ is bounded. For any critical and unbounded distribution $p$, we prove in Theorem \ref{MainA} the convergence of the GW tree with offspring distribution $p$ conditioned to have large maximal out-degree, to Kesten's tree. For any sub-critical and unbounded distribution $p$, we prove in Theorem \ref{MainB} the convergence of the conditioned GW tree with offspring distribution $p$, to condensation tree. Note that our results are complete: For a bounded distribution $p$, it is impossible to condition the GW tree to have very large maximal out-degree, see Remark \ref{bounded}; For a super-critical distribution $p$, the situation is essentially trivial, see Remark \ref{super}. Also in Theorem \ref{MainC} we study tail of the maximal out-degree of sub-critical GW trees and prove that tail of the offspring distribution $p$ and tail of the maximal out-degree are asymptotically equivalent, apart from the ratio $1-\mu_p$, where $\mu_p$ is the expectation of $p$. This result is essential for the proof of the local convergence in the sub-critical case.

On the technical level, our proofs are extremely short and elementary, thanks in particular to the convenient framework in \cite{AD14a,AD14b}. Nevertheless, we still would like to stress the following points: First, our new conditioning is the most natural way to get condensation tree in the limit, and this seems to be an intrinsic reason behind our short and elementary proofs; Second, potentially our new conditioning and some possible variants might be useful for studying condensation phenomenon in other settings;
Finally, conditionings of GW trees seem more complete now, with conditioning on large height as one extreme case, our new conditioning as the opposite extreme case, and other conditionings in between.

  This paper is organized as follows. In Section \ref{sec:Prel}, we recall several notations of trees, local convergence of random trees, and, characterizations of Kesten's tree and condensation tree. In Section \ref{sec:tail}, we study tail of the maximal out-degree of sub-critical GW trees.
  Finally Section \ref{sec:main} is devoted to the proof of our main results, Theorem \ref{MainA} for the critical case and Theorem \ref{MainB} for the sub-critical case.

\section{Preliminaries}
\label{sec:Prel}

This section is extracted from \cite{AD14a,AD14b}. For more details and proofs, refer to \cite{AD14a,AD14b}. We denote by $\mN = \{0,1,2,\ldots\}$ the set of non-negative integers and by $\mN^* = \{1,2,\ldots\}$
the set of positive integers.

\subsection{Notations of discrete trees}
\label{ss:nota}

We set
$$
\mathcal{U}=\bigcup_{n\geq 0}(\mN^*)^n
$$
the set of finite sequences of positive integers with the convention $(\mN^*)^0=\{\emptyset\}$.
If $u$ and $v$ are two sequences of $\mathcal{U}$, we denote by $uv$ the
concatenation of the two sequences, with the convention that $uv = u$ if $v = \emptyset$  and $uv = v$ if
$u = \emptyset$.
The set of ancestors of $u$ is the set:
$$
A_u = \{v \in \mathcal{U}; \mbox{ there exists }w \in \mathcal{U}, w \neq \emptyset, \mbox{ such that }u = vw\}.
$$

A tree $\t$ is a subset of $\cal{U}$ that satisfies:
\begin{itemize}
 \item[$\bullet$] $\emptyset \in \t$.

 \item[$\bullet$] If $u\in \t$, then  $A_u\subset \t$.

 \item[$\bullet$] For every $u\in \t$, there exists $k_u(\t) \in \mN \cup \{+\infty\}$ such that, for every positive integer
$i, ui \in \t$ iff $1 \leq i \leq k_u(\t)$.
 \end{itemize}

The vertex $\emptyset$ is called the root of $\t$. The integer $k_u(\t)$ represents the number of offsprings of the vertex $u\in \t$, and we call it the out-degree of the vertex $u$.
The vertex $u \in \t$ is called a leaf if $k_u(\t) = 0$ and it is said infinite if $k_u(\t) = +\infty$.
The maximal out-degree $M(\t)$ of a tree $\t$ is defined by
\beqlb\label{def:max}
M(\t)=\sup\{k_u(\t);u\in \t\}.
\eeqlb
The set of leaves of the tree $\t$ is $\mathcal{L}_0(\t) = \{u \in \t; k_u(\t) = 0\}$.

For $u \in \t$, we denote the sub-tree of $\t$ ``above'' $u$ by $\cal{S}_u(\t)$.
For $u \in \t\setminus \cal{L}_0(\t)$, denote the forest $\cal{F}_u(\t)$ ``above'' $u$ by $\cal{F}_u(\t)$.
Note that $u\in \cal{S}_u(\t)$ and $u\notin \cal{F}_u(\t)$.
For $u \in \t\setminus \{\emptyset\}$, we also define the sub-tree $\cal{S}^u(\t)$ of $\t$ ``below'' $u$ as:
$$
\cal{S}^u(\t) = \{v \in \t, u \notin A_v\}.
$$
We denote by $\mT_\infty$ the set of trees, by $\mT$ the subset of trees with no infinite vertex,
by $\mT_0$ the subset of finite trees,
by $\mT_1$ the subset of trees with a unique infinite spine but no infinite vertex,
and by $\mT_2$ the subset of trees with a unique infinite vertex but no infinite spine.
Following the terminology of Galton-Watson processes, we also say that $\t$ is extinct if $\t\in\mT_0$, non-extinct otherwise.
For precise definitions of all these notations in the present paragraph, refer to Section 2 in \cite{AD14a} and Section 2 in \cite{AD14b}.

\subsection{Local convergence of random trees}
\label{ss:lc}
For the general framework of local convergence of trees, refer to Section 2 in \cite{AD14a} and Section 2 in \cite{AD14b}.

If $\t, \s \in \mT$ and $x \in \cal{L}_0(\t)$ we denote by:
\beqlb\label{def:graft}
\t \circledast (\s,x) = \{u\in \t\}\cup \{xv, v\in \s\}
\eeqlb
the tree obtained by grafting the tree $\mathbf{s}$ on the leaf $x$ of the tree $\mathbf{t}$. For every $\mathbf{t} \in \mT$ and
every $x \in \mathcal{L}_0(\mathbf{t})$, we shall consider the set of trees obtained by grafting a tree on the
leaf $x$ of $\mathbf{t}$:
\beqlb\label{def:set}
\mT(\mathbf{t}, x) = \{\mathbf{t} \circledast (\mathbf{s},x),\mathbf{s}\in \mT\}.
\eeqlb
We recall Lemma 2.1 in \cite{AD14a}, which is a very useful characterization of convergence
in distribution in $\mT_0\cup \mT_1$.

\begin{lemma}\label{lemA}
Let $(T_n, n \in \mN^*)$ and $T$ be $\mT$-valued random variables which belong a.s.
to $\mT_0\cup \mT_1$. The sequence $(T_n, n \in \mN^*)$ converges in distribution to $T$ if and only if for
every $\mathbf{t} \in \mT_0$ and every $x \in \mathcal{L}_0(\mathbf{t})$, we have:
$$
\lim_{n\rightarrow +\infty}\mP(T_n \in \mT(\mathbf{t},x)) = \mP(T \in \mT(\mathbf{t},x))\quad and \quad
\lim_{n\rightarrow +\infty}\mP(T_n =\mathbf{t}) = \mP(T =\mathbf{t}).
$$
\end{lemma}

If $v =(v_1,\ldots, v_n) \in \cal{U}$, with $n > 0$, and $k \in\mN$, we define the shift of $v$ by $k$ as
$\theta(v, k) = (v_1+k, v_2,\ldots, v_n)$. If $\t \in \mT_0$, $\mathbf{s} \in \mT_\infty$ and $x \in \t$ we denote by:
\beqlb\label{def:graft'}
\t \circledast (\s,x) = \t \cup \{x\theta(v, k_x(\t)), v\in \s\backslash\{\emptyset\}\}
\eeqlb
the tree obtained by grafting the tree $\s$ at $x$ on ``the right'' of the tree $\t$, with the convention
that $\t\circledast(\s, x) = \t$ if $\s = \{\emptyset\}$ is the tree reduced to its root.
For every $\t\in  \mT_0$ and every $x\in  \t$, we consider the set of trees obtained by grafting a tree
at $x$ on ``the right'' of $\t$:
$$
\mT(\t, x) = \{\t \circledast (\s, x), \s \in \mT_\infty\}
$$
as well as for $k \in \mN$:
\beqlb\label{def:set'}
\mT_+(\t, x, k) = \{\s \in \mT(\t, x); k_x(\s) \geq k\}
\eeqlb
the subsets of $\mT(\t, x)$ such that $x$ has at least $k$ offsprings.
We recall Lemma 2.2 in \cite{AD14b}, which is a very useful characterization of convergence
in distribution in $\mT_0\cup \mT_2$.

\begin{lemma}\label{lemB}
Let $(T_n, n \in \mN^*)$ and $T$ be $\mT_\infty$-valued random variables which belong a.s.
to $\mT_0\cup \mT_2$. The sequence $(T_n, n \in \mN^*)$ converges in distribution to $T$ if and only if for
every $\mathbf{t} \in \mT_0$ and every $x \in \t$ and $k\in \mN$ , we have:
$$
\lim_{n\rightarrow +\infty}\mP(T_n \in \mT_+(\mathbf{t},x,k)) = \mP(T \in \mT_+(\mathbf{t},x,k))\quad and \quad
\lim_{n\rightarrow +\infty}\mP(T_n =\mathbf{t}) = \mP(T =\mathbf{t}).
$$
\end{lemma}

\subsection{GW trees}
\label{ss:GW}
Let $p=(p_0,p_1,p_2,\ldots)$ be a probability distribution on the set of the nonnegative
integers. We denote by $\mu_p$ the expectation of $p$ and assume that $0<\mu_p<+\infty$.

A $\mT$-valued random variable $\tau$ is a Galton-Watson (GW) tree with offspring distribution
$p$ if the distribution of $k_\emptyset(\tau)$ is $p$ and for $n\in \mN^*$, conditionally on $\{k_\emptyset(\tau)=n\}$, the sub-trees
$(\cal {S}_1(\tau), \cal {S}_2(\tau),\ldots, \cal {S}_n(\tau))$ are independent and distributed as the original tree $\tau$.
In particular, the restriction of the distribution of $\tau$ on the set $\mT_0$ is given by:
\beqlb\label{multi}
\forall \t \in \mT_0,\quad \mP(\tau = \t) =\prod_{u\in \t} p(k_u(\t)).
\eeqlb
The GW tree is called critical (resp. sub-critical, super-critical) if $\mu_p = 1$ (resp. $\mu_p < 1$,
$\mu_p >1$). In the critical and sub-critical case, we have that a.s. $\tau$ belongs to $\mT_0$.

Recall (\ref{def:max}), the maximal out-degree $M(\t)$ of a tree $\t$. Consider the random variable $M(\tau)$ of a GW tree $\tau$. We use $q=(q_0,q_1,q_2,\ldots)$ to denote the probability distribution of $M(\tau)$,
$H(n)$ the distribution function of $q$, and $\bar{H}(n)=1-H(n)$ the tail function.
Let $\mP_k$ be the distribution of the forest $\tau^{(k)} = (\tau_1,\ldots,\tau_k)$ of i.i.d. GW trees with offspring
distribution $p$. We set:
\beqlb\label{def:forest}
M(\tau^{(k)})=\sup_{1\leq j \leq k} M(\tau_j).
\eeqlb
When there is no confusion, we shall write $\tau$ for $\tau^{(k)}$, and $M(\tau)$ for $M(\tau^{(k)})$.

\subsection{Kesten's tree and condensation tree}
\label{ss:Kc}
Kesten's tree can be defined for critical or sub-critical offspring distributions, while condensation tree for sub-critical distributions. However in this paper we only need Kesten's tree in the critical case and condensation tree in the sub-critical case, which can be defined in a unified way.
We recall the following unified definition from Section 1 in \cite{AD14b}, which first appeared in Section 5 of \cite{J12}. Let $p$ be a critical or sub-critical offspring distribution.
Let $\tau^*(p)$ denote the random
tree which is defined by:

\begin{itemize}
 \item[i)] There are two types of nodes: $normal$ and $special$.

 \item[ii)] The root is special.

 \item[iii)] Normal nodes have offspring distribution $p$.

 \item[iv)] Special nodes have offspring distribution the biased distribution $\tilde{p}$ on $\mN\cup \{+\infty\}$
defined by:
\begin{equation*}
\tilde{p}_k =
\begin{cases}
kp_k & \text{if } k \in \mN,\\
1-\mu_p & \text{if } k = +\infty.
\end{cases}
\end{equation*}
  \item[v)] The offsprings of all the nodes are independent of each others.

   \item[vi)] All the children of a normal node are normal.

    \item[vii)] When a special node gets a finite number of children, one of them is selected uniformly
at random and is special while the others are normal.

     \item[viii)] When a special node gets an infinite number of children, all of them are normal.
 \end{itemize}
Notice that:
\begin{itemize}
 \item[$\bullet$] If $p$ is critical, then a.s. $\tau^*(p)$ has exactly one infinite spine and all its nodes have finite
out-degrees. We call it $\emph{Kesten's tree}$.
 \item[$\bullet$] If $p$ is sub-critical, then a.s. $\tau^*(p)$ has exactly one node of infinite out-degree and no infinite
spine. We call it $\emph{condensation tree}$.
 \end{itemize}

By (\ref{multi}) and the definition of Kesten's tree given above, one can prove
the following characterization of Kesten's tree, which is (2.8) in \cite{AD14a}:

\begin{lemma}\label{lemC}
Suppose that p is critical. The distribution of $\tau^*(p)$ is also
characterized by: a.s. $\tau^*(p)\in \mT_1$ and for $\t \in  \mT_0,\ x \in \cal{L}_0(\t)$,
$$
\mP(\tau^*(p) \in  \mT(\t,x)) = \frac{1}{p_0}\mP (\tau=\t).
$$
\end{lemma}

For $\t \in \mT_0$, $x \in \t$, we set:
$$
D(\t, x) =\frac{\mP(\tau = \cal {S}^x(\t))}{p_0}\mP_{k_x(\t)}(\tau = \cal {F}_x(\t)).
$$
For $z\in \mR$, we set $z_+ = \max(z, 0)$. Let $X$ be a random variable with distribution $p$.
By (\ref{multi}) and the definition of condensation tree given above,
one can prove the following characterization of condensation tree, which is Lemma 3.1 in \cite{AD14b}:

\begin{lemma}\label{lemD}
Suppose that p is sub-critical. The distribution of $\tau^*(p)$ is also
characterized by: a.s. $\tau^*(p)\in \mT_2$ and for $\t \in  \mT_0,\ x \in \t,\ k \in \mN$,
$$
\mP(\tau^*(p) \in  \mT_+(\t,x,k)) = D(\t, x)\l(1-\mu_p+\mE \l[(X - k_x(\t))_+\mathbf{1}_{\{X\geq k\}}\r]\r).
$$
\end{lemma}

\section{Tail of the maximal out-degree}
\label{sec:tail}

We begin with a simple absolute continuity result of $q_n$ and $p_n$. Recall that $q=(q_0,q_1,q_2,\ldots)$ is the probability distribution of $M(\tau)$.

\begin{lem} \label{ac} If $p_0>0$,
then for any nonnegative integer $n$, we have $q_n>0$ if and only if $p_n>0$.
If $p_0=0$, then $\mathbb{P}\l[M(\tau) =\sup p \r]=1$, where $\sup p=\sup \{i;p_i>0\}$.
\end{lem}

\proof First of all if $p_n=0$, of course $q_n=0$.
If $p_0>0$, then $\mathbb{P}\l[M(\tau) =0 \r]=p_0>0$.
If $p_0>0$ and $p_n>0$ for $n\geq 1$, then we see that $\mathbb{P}\l[M(\tau) =n \r]\geq p_n (p_0)^n>0$,
where $p_n (p_0)^n$ is the probability of the tree such that the root has $n$ offsprings and these $n$ offsprings are all leaves.
For the last statement, it holds trivially if $p_1=1$.
Clearly any other $p$ with $p_0=0$ is super-critical, then trivially $M=\sup p$ on the event of non-extinction.
This implies the last statement since when $p_0=0$, the GW tree $\tau$ is non-extinct a.s.
\qed

We also need a well-known result on $\bar{F}(n)$. Recall that we use $F(n)$ to denote the distribution function of $p$, and $\bar{F}(n)=1-F(n)$ the tail function.

\begin{lemma}\label{lemsmall}
For any offspring distribution $p=(p_0,p_1,p_2,\ldots)$ with finite expectation, we have
$$\bar{F}(n) = o\,(1/n).$$
\end{lemma}

Now we are ready to prove the main result of this section. Recall the notations $q_n$ and $H(n)$ introduced in Section \ref{ss:GW}.

\begin{theo}\label{MainC}
Suppose that the offspring distribution $p=(p_0,p_1,p_2,\ldots)$ is sub-critical and unbounded, then $\lim_{n\rightarrow\infty}H^n(n)=1$ and
$$
\lim_{n\rightarrow\infty}\frac{\bar{F}(n)}{\bar{H}(n)}=\lim_{n\rightarrow\infty}\frac{p_n}{q_n}=1-\mu_p,
$$
where the second limit is understood along the infinite subsequence $\{n; p_n > 0\}$.
\end{theo}

\proof
By considering the out-degree of the root, we get
\beqlb\label{basic}
H(n)=\sum_{m\leq n}p_m H^m(n).
\eeqlb
So $H(n)$ is a solution of the equation
$$
f_n(y)=\sum_{m\leq n}p_m y^m-y=0.
$$
Obviously this equation has a unique solution $y$ on $[0,1]$. Specifically, just check that
$f_n(0)=p_0>0$, $f_n(1)=\sum_{m\leq n}p_m -1 <0$, and $f'_n(y)\leq\sum_{1\leq m\leq n}m p_m - 1 < 0$ for $y$ on $[0,1]$.
By comparing (\ref{basic}) at the values $n-1$ and $n$, we see
\beqnn
q_n \ar=\ar \sum_{1\leq m\leq n-1} p_m [H^m(n)-H^m(n-1)] +p_n H^n(n) \cr
\ar=\ar \sum_{1\leq m\leq n-1} q_n p_m \sum_{1\leq i\leq m}[H^{m-i}(n)H^{i-1}(n-1)] +p_n H^n(n), \cr
\eeqnn
Or equivalently,
\beqlb\label{3+4}
q_n\l(1-\sum_{1\leq m\leq n-1} p_m \sum_{1\leq i\leq m}[H^{m-i}(n)H^{i-1}(n-1)]\r)=p_n H^n(n).
\eeqlb

From (\ref{3+4}) we get an obvious inequality
$$
q_n(1-\mu_p)\leq p_n.
$$
Combined with Lemma \ref{lemsmall}, we see $\bar{H}(n) = o\, (1/n)$, and consequently
\beqlb\label{3}
\lim_{n\rightarrow\infty}H^n(n)=1.
\eeqlb
From the obvious fact that for fixed $m$,
$$
\sum_{1\leq i\leq m}H^{m-i}(n)H^{i-1}(n-1)\leq m \quad \text{and} \quad \lim_{n\rightarrow\infty}\sum_{1\leq i\leq m}H^{m-i}(n)H^{i-1}(n-1)=m,
$$
we get
\beqlb\label{4}
\lim_{n\rightarrow\infty}\l(1-\sum_{1\leq m\leq n-1} p_m \sum_{1\leq i\leq m}[H^{m-i}(n)H^{i-1}(n-1)]\r)=1-\mu_p.
\eeqlb

For a sub-critical offspring distribution $p$, clearly $p_0>0$, so that by Lemma \ref{ac} we have $q_n>0$ if and only if $p_n>0$.
By combining (\ref{3+4}), (\ref{3}), and (\ref{4}), we get that
$$
\lim_{n\rightarrow\infty}\frac{\bar{F}(n)}{\bar{H}(n)}=\lim_{n\rightarrow\infty}\frac{p_n}{q_n}=1-\mu_p,
$$
where the second limit is understood along the infinite subsequence $\{n; p_n > 0\}$, and the first equality is automatic from the second one.
\qed

We end this section with two remarks related to Theorem \ref{MainC}.

\begin{remark}
In \cite{B11,B13}, Bertoin studied scaling limits and tails of the maximal out-degree of critical GW trees with regularly varying offspring tails. In Section 1 of \cite{B11} (the first paragraph on page 577) the maximal out-degree of sub-critical GW trees was mentioned, and combined with a classical result due to Gnedenko (see e.g. Proposition 1.11 in \cite{R87}), it is easy to see \cite{B11} contains the tail part ($\bar{F}(n)/\bar{H}(n)$) of our Theorem \ref{MainC} in the case of regularly varying offspring tails. Theorem \ref{MainC} deals with all sub-critical and unbounded offspring distributions and gives
both the ``density" part ($p_n/q_n$) and the tail part. Note that it is impossible to achieve this generality by using extreme value theory, since some offspring distributions do not belong to the domain of attraction of any
extreme value distribution, e.g. the geometric distribution (see e.g. \cite{A80}).
\end{remark}

\begin{remark}\label{compare}
It is interesting to compare our Theorem \ref{MainC} with Corollary 1 in \cite{B11} and Lemma 1 in \cite{B13}. In the critical and infinite variance case considered in Corollary 1 of \cite{B11}, $\bar{H}(n)$ is asymptotically larger than $\bar{F}(n)$, and $M(\tau)$ has infinite expectation. In the critical and finite variance case considered in Lemma 1 of \cite{B13}, $\bar{H}(n)$ is still asymptotically larger than $\bar{F}(n)$, and $M(\tau)$ may have finite or infinite expectation, depending on $\bar{F}(n)$. In the sub-critical case, our Theorem \ref{MainC} shows that the two tails are asymptotically of the same order, and $M(\tau)$ has finite expectation.
\end{remark}

\section{Main results}
\label{sec:main}

If a critical or sub-critical offspring distribution $p$ is unbounded, then clearly $p_0>0$, so that by Lemma \ref{ac} we have $q_n>0$ if and only if $p_n>0$.

\begin{theo}\label{MainA}
Suppose that the offspring distribution $p=(p_0,p_1,p_2,\ldots)$ is critical and unbounded,
then for the GW tree $\tau$ with the offspring distribution $p$, we have as $n\rightarrow\infty$,
$$
\normalfont\text {dist }(\tau\big|M(\tau)=n)\rightarrow\text{dist }(\tau^*(p)),
$$
where the limit is understood along the infinite subsequence $\{n; p_n > 0\}$, and as $n\rightarrow\infty$,
$$
\normalfont\text{dist }(\tau\big|M(\tau)>n)\rightarrow\text{dist }(\tau^*(p)).
$$
\end{theo}

\proof Since $p$ is critical, we have a.s. $\tau \in \mT_0$ and $\tau^*(p)\in \mT_1$.
So we will use Lemma \ref{lemA} to prove the convergence.

Recall (\ref{def:graft}) and (\ref{def:set}), the definitions of $\t \circledast (\tilde{\t},x)$ and $\mT(\t,x)$.
Using (\ref{multi}), we have for every $\t \in \mT_0$, $x \in \cal{L}_0(\t)$, and $\tilde{\t} \in \mT_0$:
$$
\mP(\tau=\t \circledast (\tilde{\t},x)) = \frac{1}{p_0}\mP (\tau=\t)\mP (\tau=\tilde{\t}).
$$
Recall (\ref{def:max}), the definition of the maximal out-degree. For the tree $\t \circledast (\tilde{\t},x)$, when $n>M(\t)$ it is clear that $M(\t \circledast (\tilde{\t},x))=n$ if and only if $M(\tilde{\t})=n$. So for such $n$ we get:
\beqnn
\mP(\tau \in \mT(\t,x),M(\tau)=n) \ar=\ar \sum_{\tilde{\t}\in\mT_0} \mP(\tau=\t \circledast (\tilde{\t},x))\mathbf{1}_{\{M(\t \circledast (\tilde{\t},x))=n\}} \cr
\ar=\ar \sum_{\tilde{\t}\in\mT_0} \frac{1}{p_0}\mP(\tau=\t) \mP (\tau=\tilde{\t})\mathbf{1}_{\{M(\tilde{\t})=n\}}\cr
\ar=\ar \frac{1}{p_0}\mP(\tau=\t)\sum_{\tilde{\t}\in\mT_0}  \mP (\tau=\tilde{\t})\mathbf{1}_{\{M(\tilde{\t})=n\}}\cr
\ar=\ar \mP(\tau^*(p) \in \mT(\t,x)) \mP (M(\tau)=n),\cr
\eeqnn
where we used Lemma \ref{lemC} for the last equality. Therefore we have
\beqnn
\mP(\tau \in \mT(t,x)\big|M(\tau)=n) \ar=\ar \frac{\mP(\tau \in \mT(\t,x),M(\tau)=n)}{\mP (M(\tau)=n)} \cr
\ar=\ar \frac{\mP(\tau^*(p) \in \mT(\t,x))\mP (M(\tau)=n)}{\mP (M(\tau)=n)} \cr
\ar=\ar \mP(\tau^*(p) \in \mT(\t,x)). \cr
\eeqnn
For all $\t \in \mT_0$ and all $n > M(\t)$, we have
$$
\mP(\tau = \t, M(\tau)=n) = \mP(\tau = \t, M(\t)=n) \leq \mathbf{1}_{\{M(\t)=n\}} = 0
$$
and thus:
$$
\mP(\tau = \t\big| M(\tau)=n) = 0= \mP(\tau^*(p) = \t).
$$
By Lemma \ref{lemA}, we have proved the first convergence. Finally since dist $(\tau|M(\tau)>n)$ is a mixture of dist $(\tau|M(\tau)=k)$
for $k > n$, the second convergence follows from the first one.
\qed

We have several remarks related to Theorem \ref{MainA}.

\begin{remark}
The above proof is essentially copied from the proof of Theorem 3.1 in \cite{AD14a}. We include the complete proof here for the readers' convenience.
Recall the function $D$ on $\mT$ defined in Section 3 of \cite{AD14a}. It seems that Theorem 3.1 in \cite{AD14a} is stated for the case of positive $D$, while our conditioning corresponds to the degenerate case of $D=0$. Clearly when $D=0$ the proof of Theorem 3.1 in \cite{AD14a} holds without assumption (3.2) there. However for our result here, we do need to know the simple fact that $q_n>0$ if and only if $p_n>0$, which is implied by our Lemma \ref{ac}. Clearly for our result the assumption $q_n>0$ for $n$ large enough is not needed, as stated in Theorem 3.1 of \cite{AD14a}.
We only need the offspring distribution $p$ to be unbounded, then interpret the first limit in Theorem \ref{MainA} along the infinite subsequence $\{n; p_n > 0\}$.
\end{remark}

\begin{remark}
If we repeat the above proof for a sub-critical and unbounded offspring distribution $p$, then we can easily confirm that dist $(\tau|M(\tau)=n)$ does not converge to the distribution of Kesten's tree associated with $p$ (see Section 2.4 in \cite{AD14a} for the definition and (2.8) in \cite{AD14a} for a characterization). We will prove in Theorem \ref{MainB} that dist $(\tau|M(\tau)=n)$ converges to the distribution of condensation tree associated with $p$, which is of course different from Kesten's tree associated with $p$.
\end{remark}

\begin{remark}
\label{rem:Kes}
Recall that Kesten's tree associated with a critical offspring distribution $p$ has no infinite node. However since there is an infinite spine and we grow on this spine infinite many i.i.d. trees associated with the offspring distribution $p'=(p'_0,p'_1,p'_2,\ldots)=(p_1,2p_2,3p_3,\ldots)$ (the roots of those trees have offspring distribution $p'$), we see a.s. $M(\tau^*(p))=m$, where $m=\sup\{n;p_n>0\}$. Then clearly $m=+\infty$ and a.s. $M(\tau^*(p))=+\infty$ if and only if $p$ is unbounded.
\end{remark}

\begin{theo}\label{MainB}
Suppose that the offspring distribution $p=(p_0,p_1,p_2,\ldots)$ is sub-critical and unbounded, then for the GW tree $\tau$ with the offspring distribution $p$, we have as $n\rightarrow\infty$,
$$
\normalfont\text{dist }(\tau\big|M(\tau)=n)\rightarrow\text{dist }(\tau^*(p)),
$$
where the limit is understood along the infinite subsequence $\{n; p_n > 0\}$, and as $n\rightarrow\infty$,
$$
\normalfont\text{dist }(\tau\big|M(\tau)>n)\rightarrow\text{dist }(\tau^*(p)).
$$
\end{theo}

\proof As in the proof of Theorem \ref{MainA}, we only need to prove the first convergence. Since $p$ is sub-critical, we have a.s. $\tau \in \mT_0$ and $\tau^*(p)\in \mT_2$. So we will use Lemma \ref{lemB} to prove the convergence.

Recall (\ref{def:graft'}), (\ref{def:set'}), and (\ref{def:forest}), the definitions of $\t \circledast (\tilde{\t},x)$, $\mT_+(\mathbf{t},x,k)$, and the maximal out-degree of a forest.
For a tree $\s\in \mT_+(\t,x,k)$, when $n>M(\t)$ it is clear that $M(\mathbf{s})=n$ if and only if $k_x(\mathbf{s})=n$ and the attached forest has maximal out-degree at most $n$, or $k_x(\mathbf{s})<n$ and the attached forest has maximal out-degree exactly $n$. Then similar to the proof of Lemma \ref{lemD},
for $k\in\mN$, $\t\in \mT_0$, $x\in \t$, $\ell = k_x(\t)$, and $n>M(\t)$, we have
\beqlb \label{1+2}
\ar\ar \mP(\tau \in \mT_+(\t,x,k),M(\tau)=n) \cr
\ar\ar \quad =\sum_{\s\in \mT_+(\t,x,k)} \mP(\tau=\s)\mathbf{1}_{\{M(\s)=n\}}\cr
\ar\ar \quad = D(\t,x)\l(p_n\mP_{n-\ell}(M(\tau)\leq n)+\sum_{j\geq\max(\ell+1,k)}^{n-1} p_j\mP_{j-\ell}(M(\tau)=n)\r).
\eeqlb
Recall the notations $q_n$ and $H(n)$ introduced in Section \ref{ss:GW}. By Theorem \ref{MainC}, we have
$$
\lim_{n\rightarrow\infty} \mP_{n-\ell}(M(\tau)\leq n)=\lim_{n\rightarrow\infty}H^{n-\ell}(n)=1,
$$
and
$$
\lim_{n\rightarrow\infty}\frac{p_n}{\mP(M(\tau)= n)}=\lim_{n\rightarrow\infty}\frac{p_n}{q_n}=1-\mu_p.
$$
So
\beqlb \label{1}
\lim_{n\rightarrow\infty}\frac{p_n\mP_{n-\ell}(M(\tau)\leq n)}{\mP(M(\tau)= n)}=1-\mu_p.
\eeqlb
Since
\beqnn
\mP_m(M(\tau)=n) \ar=\ar \mP_m(M(\tau)\leq n)-\mP_m(M(\tau)\leq n-1) \cr
\ar=\ar H^m(n)-H^m(n-1) \cr
\ar=\ar q_n\sum_{1\leq i\leq m}[H^{m-i}(n)H^{i-1}(n-1)],\cr
\eeqnn
we see
$$
\frac{\mP_m(M(\tau)=n)}{\mP(M(\tau)= n)}\leq m \quad \mbox{and}\quad \lim_{n\rightarrow\infty}\frac{\mP_m(M(\tau)=n)}{\mP(M(\tau)= n)}=m.
$$
So
\beqlb \label{2}
\lim_{n\rightarrow\infty}\frac{\sum_{j\geq\max(\ell+1,k)}^{n-1} p_j\mP_{j-\ell}(M(\tau)=n)}{\mP(M(\tau)= n)}=\sum_{j\geq\max(\ell+1,k)}(j-\ell) p_j.
\eeqlb
Combining (\ref{1+2}), (\ref{1}), and (\ref{2}) together, we have that as $n\rightarrow +\infty$,
\beqnn
\mP(\tau \in \mT_+(\mathbf{t},x,k)\big|M(\tau)=n) \ar=\ar \frac{\mP(\mT_+(\mathbf{t},x,k),M(\tau)=n)}{\mP (M(\tau)=n)} \cr
\ar\rightarrow\ar D(\t,x) \l(1-\mu_p+\sum_{j\geq\max(\ell+1,k)}(j-\ell) p_j\r)\cr
\ar=\ar \mP(\tau^*(p) \in \mT_+(\mathbf{t},x,k)), \cr
\eeqnn
where we used Lemma \ref{lemD} for the last equality.
Finally same as in the proof of Theorem \ref{MainA}, for $n > M(\t)$,
$$
\mP(\tau = \t\big| M(\tau)=n) = 0= \mP(\tau^*(p) = \t).
$$
By Lemma \ref{lemB}, we are done.
\qed

We end this section with several remarks related to both Theorem \ref{MainA} and Theorem \ref{MainB}.

\begin{remark}
Since a.s. any condensation tree $\tau^*(p)$ has an infinite node, we see that a.s. $M(\tau^*(p))=+\infty$. However, even for a slightly different maximal out-degree $M'$ defined by $M'(\t)=\sup\{k_u(\t);u\in \t, k_u(\t)<+\infty\}$, we still have that a.s. $M'(\tau^*(p))=+\infty$ for any condensation tree $\tau^*(p)$ with an unbounded offspring distribution $p$. This is because on the infinite node, we have infinite many i.i.d. trees with the unbounded offspring distribution $p$. It is interesting to compare the maximal out-degree of condensation tree with that of Kesten's tree, and we note that $M'(\tau^*(p))=+\infty$ for any Kesten's tree $\tau^*(p)$ with an unbounded offspring distribution $p$, even though Kesten's tree has no infinite node, see Remark \ref{rem:Kes}.
\end{remark}

\begin{remark}\label{bounded}
Clearly $\mP(M(\tau)=n)=0$ if $p_n=0$, so for a bounded distribution $p$, it is impossible to formally define elementary conditional probabilities with respect to the null event
$\{M(\tau)=n\}$, for all large enough $n$.
Regarding the local convergence of large conditioned GW trees, first recall that
the maximal out-degree $M$ is a function from $\mT_\infty$ to $\mN\cup\{+\infty\}$. It is easy to check that $M$ is continuous at $+\infty$ with respect to the topology induced by the distance $d_\infty$ on $\mT_\infty$ and the usual topology on $\mN\cup\{+\infty\}$.
Refer to Section 2 in \cite{AD14b} for the precise definition of $d_\infty$. Combined with Theorem 4.27 in \cite{K02} and the fact that a.s. $M(\tau^*(p))=+\infty$ for any condensation tree $\tau^*(p)$, we see that any large conditioned GW tree with a bounded and sub-critical offspring distribution $p$ does not converge in distribution to condensation tree, regardless of the conditioning. In comparison, large conditioned GW tree with a bounded and critical or sub-critical offspring distribution may converge in distribution to Kesten's tree under suitable conditionings (even include some unbounded and sub-critical distributions), see Proposition 4.1, Proposition 4.2, and Corollary 5.7 in \cite{AD14a}, and Theorem 1.3 in \cite{AD14b}.
\end{remark}

\begin{remark}\label{super}
It is well-known that a super-critical GW tree $\tau$ with offspring distribution $p$ can be decomposed into two parts: the part that $\tau$ is extinct and the part that $\tau$ is non-extinct. See e.g. Section 12 in Chapter 1 of \cite{AN72}.
When conditioned on non-extinction, a.s. $M(\tau)=+\infty$ and trivially
dist $(\tau|M(\tau)>n)$ converges to the distribution of $\tau$ conditioned on non-extinction.
The GW tree $\tau$ conditioned on extinction is a sub-critical GW tree with offspring distribution $\hat{p}$ (see e.g. Theorem 3 on page 52 in \cite{AN72}), so by Theorem \ref{MainB} we have that
as $n\rightarrow\infty$,
$\text{dist }(\tau|M(\tau)=n)\rightarrow\text{dist }(\tau^*(\hat{p}))$.
\end{remark}

\bigskip
{\bf Acknowledgements.} Sincere thanks to anonymous referees for their comments and suggestions, which improved considerably the presentation of this paper.


\end{document}